\documentclass[12pt]{article}

\usepackage{ amsmath, amsthm, amssymb, amsfonts}
\newcommand\inv{^{-1}}

 \newcommand\ch{\raisebox{.4ex}{$\chi$}}

\newcommand\iy{\infty}

 \def\({\left(} \def\){\right)}

\newcommand{\ra}{\rightarrow}
 
\newcommand{\bP}{\mathbb{P}}
\newcommand{\bR}{\mathbb{R}}

\newcommand{\cA}{\mathcal{A}}
\newcommand{\cO}{\mathcal{O}}
\newcommand{\cT}{\mathcal{T}}

\newcommand{\be}{\begin{equation}}
\newcommand{\ee}{\end{equation}}
\newcommand{\ba}{\begin{eqnarray*}}
\newcommand{\ea}{\end{eqnarray*}}
\newcommand{\bae}{\begin{eqnarray}}
\newcommand{\eae}{\end{eqnarray}}
\newcommand{\bc}{\begin{center}}
\newcommand{\ec}{\end{center}}

\newcommand{\fr}{\frac}

\newcommand{\ha}{\fr{1}{2}\,}

\newcommand{\ai}{\textrm{Ai}}
\def\ch{\raisebox{.3ex}{$\chi$}}

\newcommand{\tr}{\ensuremath{\textrm{tr }}}
\newcommand{\airy}{\ensuremath{ \textrm{Ai} }}
\newcommand{\airytwo}{\ensuremath{\mathcal{A}_2}}
\newcommand{\bmb}{\ensuremath{\begin{bmatrix}}}
\newcommand{\emb}{\ensuremath{\end{bmatrix}}}

\begin{document}
\vspace*{-15ex}
\textsc{\hfill  February 8, 2011}
\begin{center}{\large \bf Asymptotics for the Covariance of the Airy$_2$ Process}\end{center}

\begin{center}{\large \bf Gregory Shinault}\\
{\it Department of Mathematics\\
University of California\\
Davis, CA 95616, USA\\
email: gshinault@math.ucdavis.edu}\end{center}

\begin{center}{\large\bf Craig A.~Tracy}\\
{\it Department of Mathematics \\
University of California\\
Davis, CA 95616, USA\\
email: tracy@math.ucdavis.edu}\end{center}

\begin{abstract}
In this paper we compute some of the higher order terms in the asymptotic behavior of the two point function $\bP(\airytwo(0)\leq s_1,\airytwo(t)\leq s_2)$, extending the previous work of Adler and van Moerbeke \cite{AvM1, AvM2} and Widom \cite{Wi}. We prove that it is possible to represent any order asymptotic approximation as a polynomial and integrals of the Painlev\'e II function $q$ and its derivative $q'$.  Further, for up to tenth order we give this asymptotic approximation as a linear combination of the Tracy-Widom GUE density function $f_2$ and its derivatives.  As a corollary to this, the asymptotic covariance is expressed up to tenth order in terms of the moments of the Tracy-Widom GUE distribution.
\end{abstract}

\bc \S1. \textsc{Introduction}\ec
\par
The \textit{Airy}$_2$ \textit{process}, $\cA_2(t)$,  introduced by Pr\"ahofer and Spohn \cite{PrSp} in the context of the 
polynuclear growth (PNG) model, is a stationary stochastic process whose joint distributions for $t_1<\cdots< t_m$ are given by
\be \bP\left(\cA_2(t_1)\le s_1,\ldots, \cA_2(t_m)\le s_m\right)=\det\left(I-\chi K_2 \chi\right)_{L^2(\{1,\ldots,m\}\times\bR)} \label{airy2}\ee
where $K_2\doteq K_2(x,y)$ is a $m\times m$ matrix  kernel, called the \textit{extended Airy kernel}, given by
\[ K_{i,j}(x,y)=\left\{ \begin{array}{rl}
\int_0^\iy e^{-z(t_i-t_j)}\, \ai(x+z)\ai(y+z)\, dz& \textrm{if}\>\> i\ge j,\\
-\int_{-\iy}^0 e^{-z(t_i-t_j)}\, \ai(x+z) \ai(y+z)\, dz & \textrm{if}\>\> i<j,
\end{array}\right.\]
and $\ch=\ch(x)$ is the $m\times m$ diagonal matrix whose $i$th diagonal term is the indicator function $\ch_{s_i}(x)=1_{x>s_i}$.
The one-point function $\bP(\cA_2(t)\le s)$ is the Tracy-Widom  GUE distribution, $F_2(s)$ \cite{TWAiry}.

For stochastic growth models and their  closely related interacting particle systems, the Airy$_2$ process is fundamental since it is expected to describe the limiting process for height fluctuations belonging to  the \textit{KPZ Universality Class} with droplet (or step) initial conditions (see
\cite{Do2, SS3} for recent reviews).  This has been proved for the PNG model \cite{PrSp}, the discrete PNG model \cite{Jo1}, the boundary of the north polar region of the Aztec diamond \cite{Jo2} and the totally asymmetric simple exclusion process (TASEP) \cite{Jo2}.  At the the level of the one-point function, this universality has been established for the asymmetric simple exclusion process (ASEP) \cite{TWAsep} and for the KPZ equation \cite{SS1, SS2, ACQ, Do1}.  Recent work \cite{ProlSp} using replica methods have extended this KPZ work to the 2-point function.  The Airy$_2$ process also describes the limiting process of the largest eigenvalue in Dyson's Brownian motion model in random matrix theory.  For further appearances of the Airy$_2$ process see \cite{BFS, BF, FS, IS}.

We summarize some known properties of the Airy$_2$ process:
\vspace{-2ex}
\begin{enumerate}
\item $\cA_2(t)$ has continuous sample paths \cite{PrSp, Jo1}.\vspace{-1ex}
\item $\cA_2(t)$ locally looks like Brownian motion (see H\"agg \cite{Ha} for a precise statement).\vspace{-1ex}
\item The distribution functions (\ref{airy2}) satisfy nonlinear differential equations \cite{AvM1, AvM2, TWDyson}. \vspace{-1ex}
\item The covariance $\textrm{cov}_2(t):=\textrm{cov}\left(\cA_2(t)\cA_2(0)\right)$ has the following asymptotic expansions:
\be \textrm{cov}_2(t)=\left\{\begin{array}{ll} \textrm{var}(F_2) - t +\cO(t^2), & t\ra 0^+, \\ 
		 \displaystyle\sum_{n=1}^N \dfrac{C_{n}}{t^{n}} +\cO\left(\dfrac{1}{t^{N+1}}\right),& t\ra\iy. \end{array}\right.\label{cov}\ee
\end{enumerate}

The small-$t$ expansion of $\textrm{cov}_2$ was given by Pr\"ahofer and Spohn \cite{PrSp} and they also found the leading large-$t$ term $C_1=0$, $C_2=1$.  The existence of the higher order terms in the large-$t$ expansion
of $\textrm{cov}_2$ was established by Adler and van Moerbeke \cite{AvM1, AvM2} and by Widom \cite{Wi} using different methods.   In both
\cite{AvM1, AvM2} and \cite{Wi} the coefficient $C_4$ was expressed as a double integral whose integrand was in terms of the Hastings-McLeod solution of Painlev\'e II appearing in the distribution $F_2$.  One of the main results of this paper is to
prove that the coefficients $C_{2n}$, $2\le n \le 5$, are expressible in terms of the moments of $F_2$.  (The odd coefficients are all zero.)
Precisely, if
\[ \mu_n:=\int_{-\iy}^\iy s^n f_2(s)\, ds,\>\>f_2:=F_2', \]
then
\bae
C_4&=&2\mu_1, \label{1stC4} \\
C_6&=& 2\mu_2 +\fr{10}{3} \mu_1^2, \label{1stC6} \\
C_8&=& 2\mu_3 + 14\mu_2\mu_1+\fr{13}{2}, \label{1stC8} \\
C_{10}&=& 2\mu_4 + 24 \mu_3 \mu_1 +\fr{126}{5}\mu_2^2+ 116 \mu_1. \label{1stC10}
\eae
We conjecture that $C_{2n}$ can be expressed in terms of a  polynomial in  $\mu_i$ for $i\le n-1$.  To prove these results we follow the program established by Widom \cite{Wi} and first prove
\bae
\bP\left( \cA_2(t) \le s_2, \cA_2(0) \le s_1 \right) = \sum_{n=0}^N \fr{c_n(s_1,s_2)}{t^n} + \cO\left( \fr{1}{t^{N+1}} \right) \label{twopoint}
\eae
as $t \to \iy$. It was previously shown \cite{AvM1, AvM2, Wi} that for $n\leq 4$ each $c_n$ could be written as polynomials and integrals of the Painlev\'e II function, its derivative, and the variables $s_1$ and $s_2$.  A feature of our analysis is that we show through order $t^{-10}$ that each $c_{2n}$ can be expressed in terms of $f_2$, its derivatives, and polynomials in $s_k$; see \eqref{c2}--\eqref{c6}, \eqref{c8}.

Bornemann \cite{Bo1, Bo2} has given a high precision numerical evaluation of $\textrm{cov}_2(t)$, $0\le t\le 100$.  His method involves a numerical evaluation of the Fredholm determinant appearing in (\ref{airy2}) for $m=2$ followed by  numerical integrations to give $\textrm{cov}_2(t)$. In Appendix 2 we compare the large-$t$ asymptotics with these numerical results. 

In the present paper we begin by showing how to obtain an asymptotic expression for  the extended Airy kernel following Widom \cite{Wi}.  The large-$t$ asymptotics of the two-point distribution is then given in terms of $f_2$ and its derivatives.  This in turn allows the easy computation of the large-$t$ expansion of the covariance, which ends the main body of the paper.  The appendices contain some of the higher order terms, and comparison to high precision numerical results.

\bc \S2. \textsc{Asymptotics for $\ch K_2 \ch$}\ec
\par
The first step in our asymptotic analysis is a large $t$ expression for the extended Airy kernel $\ch K_2 \ch$.  In the $m=2$ case the $\ch K_2 \ch$ operator has a matrix kernel of the form
\[
\left[\begin{array}{cc}
 \ch_{s_1}(x) \displaystyle\int_0^\infty \airy(x+z)\airy(y+z) \ dz \ \ch_{s_1}(y) & \hspace{-.15in} \ch_{s_1}(x) \displaystyle\int_0^\infty e^{-zt} \airy(x+z)\airy(y+z) \ dz \ \ch_{s_2}(y) \\&\\
-\ch_{s_2}(x) \displaystyle\int^0_{-\infty}e^{zt} \airy(x+z)\airy(y+z) \ dz \ \ch_{s_1}(y) & \hspace{-.15in} \ch_{s_2}(x) \displaystyle\int_0^\infty \airy(x+z)\airy(y+z) \ dz \ \ch_{s_2}(y) \end{array} \right].
\]
We must compute the Fredholm determinant of this operator for large-$t$.  To this end we will split the operator into two manageable components so that $\ch K_2 \ch = K + L$:
\[
K(x,y):=\left[\begin{array}{cc}
 \ch_{s_1}(x) \displaystyle\int_0^\infty \airy(x+z)\airy(y+z) \ dz \ \ch_{s_1}(y) & 0 \\&\\
0 & \hspace{-15ex}\ch_{s_2}(x) \displaystyle\int_0^\infty \airy(x+z)\airy(y+z) \ dz \ \ch_{s_2}(y) \end{array} \right],
\]
\vspace{.05in}
\[
L(x,y) := \left[\begin{array}{cc}
 0 & \hspace{-15ex}\ch_{s_1}(x) \displaystyle\int_0^\infty e^{-zt} \airy(x+z)\airy(y+z) \ dz \ \ch_{s_2}(y) \\&\\
-\ch_{s_2}(x) \displaystyle\int^0_{-\infty}e^{zt} \airy(x+z)\airy(y+z) \ dz \ \ch_{s_1}(y) & 0 \end{array} \right] .
\]
The determinant computation is simplified via 
\begin{align*}
\det(I-\ch K_2 \ch) &= \det(I-K-L) = \det\left[(I-K)(I-(I-K)^{-1}L)\right] \\
&= \det(I-K)\det(I-(I-K)^{-1}L).
\end{align*}
The determinant of $I-K$ is $F_2(s_1)F_2(s_2)$, and has no dependence on $t$.  So we need only look at $L$ to determine the asymptotics.  For that determinant, we make an expansion of the terms in $L(x,y)$.  By repeatedly applying integration by parts, the upper-right corner is
\begin{align*}
L_{12} &= \ch_{s_1}(x)\int_0^\infty dz \ e^{-zt}\airy(x+z)\airy(y+z)\ch_{s_2}(y) \\
&= \sum_{n=0}^N  \frac{1}{t^{n+1}} \sum_{k=0}^n  {n \choose k} \ch_{s_1}\airy^{(k)}\otimes\airy^{(n-k)}\ch_{s_2} + \mathcal{O} \left(t^{-(N+1)}\right)
\end{align*}
and the lower-left corner
\begin{align*}
L_{21} &= \ch_{s_2}(x)\int_{-\infty}^0 dz \ e^{zt}\airy(x+z) \airy(y+z) \ch_{s_1} (y) \\
&= \sum_{n=0}^N  \frac{(-1)^{n+1}}{t^{n+1}} \sum_{k=0}^n  {n \choose k} \ch_{s_2}\airy^{(k)}\otimes\airy^{(n-k)}\ch_{s_1} + \mathcal{O}\left(t^{-(N+1)}\right) . 
\end{align*}
We take this approximation up to $N=10$, which is valid in the trace norm, because this is the highest order term we aim to calculate.  

Now we analyze $T:=(I-K)^{-1}L$. At this point we introduce the notation
\[ Q_{n,k} = (I-\ch_{s_k} K_{\airy}\ch_{s_k})^{-1}\ch_{s_k}\airy^{(n)}. \]
The kernel $K_{\airy}$ is the classical Airy kernel of random matrix theory, while the purpose of the $\ch$'s are to ensure the kernels are integrated  over the appropriate domain. Using this notation in the asymptotic expansion of $L(x,y)$, the upper right term in the matrix kernel of $T$ is written
\[ T_{12} = \sum_{n=0}^N  \frac{1}{t^{n+1}} \sum_{k=0}^n  {n \choose k} Q_{k,1}\otimes\airy^{(n-k)}\ch_{s_2} + \mathcal{O}(t^{-(N+1)}) \]
and the lower left term is
\begin{align*}
T_{21} &= \sum_{n=0}^N  \frac{(-1)^{n+1}}{t^{n+1}} \sum_{k=0}^n  {n \choose k} Q_{k,2}\otimes\airy^{(n-k)}\ch_{s_1} + \mathcal{O}(t^{-(N+1)}).
\end{align*}
The determinant we are now interested in is best analyzed in terms of the following trace formula:
\begin{align}
\det(I-T) &= \det(I-T_{12}T_{21}) = \exp\left[\textrm{tr} \log(I-\cT)\right] \label{det} \\
&= 1-\textrm{tr}\, \cT+{1\over 2}\Big((\textrm{tr}\,\cT)^2-\textrm{tr}\,\cT^2\Big) -{1\over6}\left( (\tr \cT)^3-3 \tr \cT^2\tr \cT + 2 \tr \cT^3 \right) +\cdots. \nonumber
\end{align}
Here we have introduced $\cT = T_{12}T_{21}$. To evaluate an expression for $\cT$ we will repeatedly use the fact that $(f\otimes g)(e\otimes h)=(g,e)_{L^2}f\otimes h$.  By introducing the notation
\[ u_{n,j}(s_k) = (Q_{n,k}, \ch_{s_k}\airy^{(j)})_{L^2} = (Q_{j,k}, \ch_{s_k}\airy^{(n)})_{L^2} = u_{j,n}(s_k)
 \] 
we are able to write
\begin{align} 
\cT &= \sum_{n_1+n_2 \leq N-2} \frac{(-1)^{n_2+1}}{t^{n_1+n_2+2}}\sum_{k_1=0}^{n_1} \sum_{k_2=0}^{n_2} {n_1 \choose k_1}{n_2 \choose k_2} u_{n_1-k_1,k_2}(s_2)Q_{k_1,1}\otimes \airy^{(n_2-k_2)}\ch_{s_1} + \mathcal{O}\left(\frac{1}{t^{N+1}}\right). \label{scriptt} 
\end{align}
From \eqref{scriptt} above, we know that we will also need the formula for $(T_{12}T_{21})^2$:
\begin{align}
\cT^2 = &\sum_{\substack{n_1+n_2\\+m_1 +m_2\leq N-4}} \frac{(-1)^{n_2+m_2}}{t^{n_1+n_2+m_1+m_2+4}} \sum_{k_1=0}^{n_1} \sum_{k_2=0}^{n_2} \sum_{\ell_1=0}^{m_1} \sum_{\ell_2=0}^{m_2} {n_1 \choose k_1}{n_2 \choose k_2}{m_1 \choose \ell_1}{m_2 \choose \ell_2} \nonumber \\
&\times u_{n_1-k_1,k_2}(s_2)u_{m_1-\ell_1,\ell_2}(s_2)u_{n_2-k_2,\ell_1}(s_1)Q_{k_1,1}\otimes \ch_{s_1} \airy^{(m_2-\ell_2)} + \mathcal{O}\left(\frac{1}{t^{N+1}}\right). \label{scriptt2}
\end{align}
We continue multiplying the sums $\cT^n$ in this manner as necessary for the order of our desired approximation.

\bc \S3. \textsc{The Two-Point Distribution}\ec
\par
From \eqref{det} we conclude $\tr \cT$ is order $t^{-2}$, while \eqref{scriptt} and \eqref{scriptt2} tell us $(\tr \cT)^2- \tr \cT^2$ is order $t^{-4}$.  Terms can be added as necessary to get any order approximation one wishes for.  In other words,
\[ \bP\left( \cA_2(t) \le s_2, \cA_2(0) \le s_1 \right) = \sum_{n=0}^N \fr{c_n(s_1,s_2)}{t^n} + \cO\left( \fr{1}{t^{N+1}} \right) \]
for $t \to \iy$.

It has already been established that $c_1=c_3=0$, and $c_2$ and $c_4$ can be written in terms of the Hastings-McLeod solution, which we denote by $q(s)$, to the Painlev\'e II equation \cite{AvM1, AvM2, Wi}. This result can be extended to any order.  It is clear from \eqref{scriptt} and \eqref{scriptt2} that all of our coefficients $c_n(s_1,s_2)$ are polynomials in the variables $u_{k,j}$ multiplied by $F_2(s_1)F_2(s_2)$.   In our notation the first two terms are
\bae
c_2(s_1,s_2) &=& F_2(s_1)u_{0,0}(s_1)F_2(s_2)u_{0,0}(s_2), \label{c2a} \\
c_4(s_1,s_2) &=& F_2(s_1)F_2(s_2)\left[ u_{1,0}(s_1)u_{1,0}(s_2) + u_{0,0}(s_1)\left(2u_{2,0}(s_2)-u_{1,1}(s_2) \right) \right] \label{c4a} \\
& \ & + \ \text{ reversed},  \nonumber
\eae
where ``reversed'' denotes the interchange of $s_1$ and $s_2$ in the previous terms.  The formula for $F_2$ in terms of $q$, 
\[ F_2(s) = \exp\left[ - \int_s^\infty (x-s)q(x)^2 \ dx \right] , \]
 is well known so we need only focus on the $u_{k,j}$'s. There are three facts that will enable us to show each $u_{k,j}$ can be expressed in terms of $q$ and $q'$ and integrals of $q$ and $q'$:
\begin{enumerate}
\item $q_0 := (I-K_{\airy})^{-1}\ch_s\airy(s)$ is a solution to Painlev\'e II ($q_0 = q$),
\item $u'_{k,j} = -q_k q_j$, where $q_j:=(I-K_{\airy})^{-1}\ch_s\airy^{(j)}(s)$,
\item $q_n = (n-2)q_{n-3} + sq_{n-2} - u_{n-2,1}\,q + u_{n-2,0}\,q_1$ for $n\geq 3$.
\end{enumerate}
The first statement is established explicitly in \cite{TWAiry, tw3}, where a template is also developed for the proof of the second statement.  We begin by proving the second statement.  Henceforth we will use the notation 
\ba
Q_k(x;s):= (I-K_{\airy})^{-1}\ch_s\airy^{(k)}. \label{Qk}
\ea

\textit{Proof of Statement 2.}  We proceed by first differentiating $u_{k,j}$.
\begin{align}
{ du_{k,j} \over ds} &= {d \over ds} \left[ \int_s^\infty Q_k(x;s)\airy^{(j)}(x) \ dx \right] = -q_k(s)\airy^{(j)}(s) + \int_s^\infty { \partial Q_k \over \partial s} \airy^{(j)}(x) \ dx. \label{diffu}
\end{align}
We now need a formula for $\partial_s Q_k$:
\begin{align*}
{ \partial Q_k \over \partial s} &= { \partial \over \partial s}(I-K_{\airy})^{-1}\ch_s \airy^{(k)} = (I-K_{\airy})^{-1}{\partial K_{\airy} \over \partial s}(I-K_{\airy})^{-1}\ch_s \airy^{(k)}.
\end{align*}
This necessitates a formula $\partial_sK_{\airy}$. Direct computation shows that this has kernel
\[ {\partial K_{\airy} \over \partial s} = -K_{\airy}(x,s)\delta(y-s) .\]
By introducing the operator $R=(I-K_{\airy})^{-1}K_{\airy}$ we can write
\[ {\partial Q_k \over \partial s}(x;s) = -R(x,s)q_k(s). \]
Returning to \eqref{diffu}, we now have
\begin{align*}
{ du_{k,j} \over ds} &= -q_k(s)\airy^{(j)}(s) - \int_s^\infty  R(s,x)q_k(s)\airy^{(j)}(x) \ dx \\
&= -q_j(s)\left( \int_s^\infty (\delta(s-x)+R(s,x))\airy^{(k)}(x) \ dx  \right) \\
&= -q_j(s)(I-K)^{-1}\ch_s\airy^{(k)}(s) = -q_j(s)q_k(s).
\end{align*}
Note that we pass from the the second line to the third line because the kernel of $(I-K_{\airy})^{-1}$ is $\delta(x-y)+R(x,y)$.  \hfill $\square$

\textit{Proof of Statement 3.}  First, we develop an explicit formula for $q_n$. We start with
\ba
Q_n(x;s) = (I-K_{\airy})\inv\ch_s\ai^{(n)}
\ea
and the formula
\ba
\ai^{(n)}(x) = (n-2)\ai^{(n-3)}(x) + x\ai^{(n-2)}(x)
\ea
for $n\geq 3$. This gives 
\ba
Q_n(x;s) = (n-2)Q_{n-3}(x;s) + P_{n-2}(x;s)
\ea
where $P_{n}(x;s):= (I-K_\ai)\inv x\ai^{(n)}(x)\ch_s$. Now we apply the commutator relations
\ba
[M,(I-K_\ai)\inv]\ai^{(n)} = xQ_n(x;s)-P_n(x;s)
\ea
and
\ba
[M,(I-K_\ai)\inv] \doteq Q_0(x;s)Q_1(y;s)-Q_1(x;s)Q_0(y;s)
\ea
to arrive at
\ba
P_n(x;s) = xQ_n(x;s)-u_{n,1}(s)Q_0(x;s) + u_{n,0}(s)Q_1(x;s).
\ea
Taking $x=s$ we get the recursive relation
\ba
q_n = (n-2)q_{n-3} + sq_{n-2} - u_{n-2,1}q + u_{n-2,0}q_1.
\ea
for $n\geq 3$. From \cite{TWFred, tw3, Wi} we know that $q_1 = q'+u_{00}\,q$ and $q_2 = sq_1 - u_{11}\,q + u_{00}\,q_1$, so we have formulas for all $q_n$. \hfill $\square$ \\

From here, we can complete the proof inductively.  If $u_{k,j}$ is expressible as a polynomial and integrals of the $q$, $q'$ and $s$ for $j\leq k$ then so is $q_{k+1}$ due to the identity above.  The differential equations for each $u_{k,j}$ give us $u_{k+1,j} = \int_s^\infty q_{k+1}(x) q_j(x) dx$ for $j\leq k+1$.  This completes the proof, so long as enough base cases are satisfied.

While it is interesting that we can write these $u_{jk}$ formulas as a polynomial and integrals in terms of $s$, $q$, and $q'$, there are explicit polynomial formulas in the aforementioned terms for each $u_{j,k}$ with $j+k\leq 8$.  For example, 
\bae
u_{0,0} &=& (q')^2 -sq^2-q^4, \label{u00} \\
u_{1,0} &=& \ha u^2_{0,0} - \ha q^2, \label{u10} \\
u_{1,1} &=& \fr{1}{3}\,u^3_{0,0}-\left(q^2+\fr{s}{3}\right)u_{0,0} - \fr{2}{3}\,q q', \label{u11} \\
u_{2,0} &=& \ha u^2_{1,1} + \ha s u_{0,0}. \label{u20}
\eae
We conjecture this is true for all $u_{j,k}$.  The conclusion is the two-point distribution has an asymptotic formula up to order $t^{-10}$  in terms of $q$.

The formulas for $q_n$ and $u'_{j,k}$ also enable us to show that $F_2u_{j,k}$ is multilinear in $f_2 = F_2'$, its derivatives, and $s$ for $j$ and $k$ needed through $c_{10}$.  A table of these formulas is included in Appendix 1.  The verification of each equation simply requires the differentiation of both sides and looking at $s\to \infty$ to verify the constants of integration are zero.

Using this table, \eqref{scriptt}, and \eqref{scriptt2}, a direct computation yields the formulas
\bae
c_2(s_1,s_2) &=& f_2(s_1)f_2(s_2), \label{c2} \\
c_4(s_1,s_2) &=& (s_1+s_2)f_2(s_1)f_2(s_2) + \frac{1}{2}f'_2(s_1)f'_2(s_2), \label{c4} \\
c_6(s_1,s_2) &=& \fr{1}{3}\,(3s_1+s_2)(3s_2+s_1)f_2(s_1)f_2(s_2)+3\left(f_2'(s_1) f_2(s_2)+f_2(s_1) f_2'(s_2)\right)
\nonumber\\
&&\hspace{2ex}+(s_1+s_2)f_2'(s_1)f_2'(s_2)+\fr{1}{6}\,f_2''(s_1) f_2''(s_2)\label{c6}.
\eae
The odd terms are not listed here because they are zero.  The coefficient formula for $c_8$ is given below in \eqref{c8}.  

From a logical standpoint the work we have done is sufficient.  However, it is instructional to see an alternate derivation of \eqref{c2}--\eqref{c4} from \eqref{c2a}--\eqref{c4a}.  First note that we have already established $u_{0,0}(s) = \int_s^\infty q(x)^2 \ dx$.  Using this and the integral formula for $F_2$, it is easy to check that $f_2=F'_2=F_2u_{0,0}$.  This gives us \eqref{c2}.  For the next term, we use the identities \eqref{u20}, \eqref{u10}, and $f_2=F_2u_{0,0}$ in the formula \eqref{c4a} to get
\[ c_4(s_1,s_2) = \ha F_2(s_1)(u_{0,0}(s_1)^2-q(s_1)^2)F_2(s_2)(u_{0,0}(s_2)^2-q(s_2)^2) + (s_1+s_2)f_2(s_1)f_2(s_2). \]
Finally, observe that $f'_2 = [F_2u_{0,0}]' = F_2u_{0,0}^2 - F_2q^2$ so we arrive at \eqref{c4}.  

\bc \S4. \textsc{The Covariance}\ec
\par
For the covariance of the Airy$_2$ process we denote the $N$-th order asymptotic approximation by
\ba
\textrm{cov}_{2,N}(t) = \sum_{n=1}^N \fr{C_{n}}{t^{n}}.
\ea
It was established by Adler and van Moerbeke \cite{AvM1, AvM2} and Widom \cite{Wi} that the coefficients up to $N=4$ are $C_1 = 0$, $C_2 = 1$, $C_3=0$, and 
\[ C_4 = \iint_{\bR^2} c_4(u,v) \ du \, dv. \]
Through an application of Fubini's theorem and \eqref{c4}, we get
\bae
C_4 &=& 2 \int_\bR uf_2(u) du \int_\bR f_2(v) dv + \frac{1}{2}\left(\int_\bR f'_2(u) du\right)^2 = 2\mu_1. \label{C4}
\eae
Following this same procedure for $C_6$ we have
\bae
C_6&=&\int_{\bR^2} c_6(s_1,s_2)\, ds_1 \, ds_2\nonumber\\
&=&\int_{\bR^2} \fr{1}{3}\,(3s_1+s_2)(3s_2+s_1)f_2(s_1)f_2(s_2)\,ds_1 \, ds_2\nonumber\\
&=& 2\mu_2 +\fr{10}{3}\mu_1^2\label{C6}
\eae
since all the other integrals appearing in \eqref{c6} integrate to zero. The analogous formulas for $C_8$ and $C_{10}$ appear in \eqref{1stC8} and \eqref{1stC10} with some details of their computation in Appendix 1.

\bc \textsc{Acknowledgements}\ec
\par
The authors acknowledge helpful conversations with Harold Widom.  We would also like to thank Folkmar Bornemann for providing us with his table of values for the Airy$_2$ covariance.  This work was supported by the National Science Foundation through grant DMS-0906387.

\pagebreak

\bc \textsc{Appendix 1: Identities}\ec
\par
In this appendix we give the formulas necessary to attain \eqref{c2}--\eqref{c6}.  The table below gives $u_{j,k}F_2$ in terms of $f_2$ and its derivatives.
\[
\begin{array}{|c|ccccccccc|}
\hline j,k  & f_2 & f_2' & f_2'' & f_2^{(3)} & f_2^{(4)} & f_2^{(5)} & f_2^{(6)} & f_2^{(7)} & f_2^{(8)} \\ \hline
0,0 & 1 & 0 & 0 & 0 & 0 & 0 & 0 & 0 & 0 \\
1,0 & 0 & \fr{1}{2} & 0 & 0 & 0 & 0 & 0 & 0 & 0 \\
1,1 & -\fr{s}{3} & 0 & \fr{1}{3} & 0 & 0 & 0 & 0 & 0 & 0 \\
2,0 & \fr{s}{3} & 0 & \fr{1}{6} & 0 & 0 & 0 & 0 & 0 & 0 \\
2,1 & -\fr{1}{4} & 0 & 0 & \fr{1}{8} & 0 & 0 & 0 & 0 & 0 \\
3,0 & \fr{7}{12} & \fr{s}{3} & 0 & \fr{1}{24} & 0 & 0 & 0 & 0 & 0 \\
2,2 & \fr{s^2}{5} & -\fr{3}{10} & 0 & 0 & \fr{1}{20} & 0 & 0 & 0 & 0 \\
3,1 & -\fr{s^2}{5} & \fr{2}{15} & \fr{s}{6} & 0 & \fr{1}{30} & 0 & 0 & 0 & 0 \\
4,0 & \fr{s^2}{5} & \fr{47}{60} & \fr{s}{6} & 0 & \fr{1}{120} & 0 & 0 & 0 & 0 \\
3,2 & \fr{2s}{9} & \fr{s^2}{18} & -\fr{1}{12} & \fr{s}{18} & 0 & \fr{1}{72} & 0 & 0 & 0 \\
4,1 & -\fr{13s}{18} & -\fr{s^2}{18} & \fr{11}{24} & \fr{s}{9} & 0 & \fr{1}{144} & 0 & 0 & 0 \\
5,0 & \fr{101s}{90} & \fr{23s^2}{90} & \fr{59}{120} & \fr{s}{18} & 0 & \fr{1}{720} & 0 & 0 & 0 \\
3,3 & -\fr{s^3}{7}+\fr{34}{63} & \fr{17s}{42} & \fr{s^2}{9} & -\fr{1}{18} & \fr{s}{36} & 0 & \fr{1}{252} & 0 & 0 \\
4,2 & \fr{s^3}{7}-\fr{11}{42} & -\fr{13s}{84} & 0 & \fr{1}{8} & \fr{s}{24} & 0 & \fr{1}{336} & 0 & 0 \\
5,1 & -\fr{s^3}{7}-\fr{74}{105} & -\fr{47s}{420} & \fr{s^2}{10} & \fr{7}{20} & \fr{s}{24} & 0 & \fr{1}{840} & 0 & 0 \\
6,0 & \fr{s^3}{7}+\fr{1151}{630} & \fr{733s}{420} & \fr{7s^2}{45} & \fr{71}{360} & \fr{s}{72} & 0 & \fr{1}{5040} & 0 & 0 \\
4,3 & -\fr{s^2}{4} & \fr{127}{288} & \fr{5s}{18} & \fr{s^2}{18} & \fr{5}{288} & \fr{s}{72} & 0 & \fr{1}{1152} & 0 \\
5,2 & \fr{43s^2}{60} & \fr{s^3}{15}-\fr{123}{160} & -\fr{3s}{20} & \fr{s^2}{30} & \fr{61}{480} & \fr{s}{60} & 0 & \fr{1}{1920} & 0 \\
6,1 & -\fr{73s^2}{60} & -\fr{s^3}{15}+\fr{271}{1440} & \fr{29s}{36} & \fr{17s^2}{180} & \fr{221}{1440} & \fr{s}{90} & 0 & \fr{1}{5760} & 0 \\
7,0 & \fr{691s^2}{420} & \fr{22s^3}{105}+\fr{4873}{1440} & \fr{394s}{315} & \fr{11s^2}{180} & \fr{83}{1440} & \fr{s}{360} & 0 & \fr{1}{40320} & 0 \\
4,4 & \fr{s^4}{9}-\fr{118s}{81} & -\fr{17s^2}{54} & -\fr{s^3}{81}+\fr{119}{144} & \fr{31s}{108} & \fr{s^2}{27} & \fr{1}{144} & \fr{s}{216} & 0 & \fr{1}{5184} \\
5,3 & -\fr{s^4}{9}+\fr{667s}{810} & \fr{197s^2}{540} & \fr{32s^3}{405}-\fr{4}{45} & \fr{119s}{1080} & \fr{29s^2}{1080} & \fr{11}{360} & \fr{11s}{2160} & 0 & \fr{1}{6480} \\
6,2 & \fr{s^4}{9}+\fr{322s}{405} & \fr{103s^2}{540} & \fr{s^3}{162}-\fr{35}{72} & \fr{113s}{540} & \fr{37s^2}{1080} & \fr{11}{180} & \fr{s}{216} & 0 & \fr{1}{12960} \\
7,1 & -\fr{s^4}{9}-\fr{21167s}{5670} & -\fr{1999s^2}{3780} & \fr{37s^3}{567}+\fr{115}{63} & \fr{6107s}{7560} & \fr{47s^2}{1080} & \fr{17}{360} & \fr{s}{432} & 0 & \fr{1}{45360} \\
8,0 & \Large \fr{s^4}{9}+\fr{19912s}{2835} & \fr{5297s^2}{1890} & \fr{409s^3}{2835}+\fr{28319}{10080} & \fr{4273s}{7560} & \fr{19s^2}{1080} & \fr{19}{1440} & \fr{s}{2160} & 0 & \fr{1}{362880} \\ \hline
\end{array}
\]
Thus, for example, $u_{2,1}F_2 = -\fr{1}{4}f_2 + \fr{1}{8}f_2^{(3)}$. To verify these identities, first divide both sides by $F_2$, differentiate both sides, then compare the behavior at infinity to check that the constant of integration is the same.  To illustrate this, we check the identity for $u_{1,0}F_2$.  After differentiation we must verify that
\[ -q_1q_0 = \fr{(f''_2/2)F_2 - (f'_2/2)f_2}{F_2^2}, \]
which is equivalent to
\[ -(q'+u_{0,0}q)qF_2 = \ha (-2qq'F_2 -2q^2f_2+u_{0,0}f'_2) - \ha f'_2 u_{0,0}. \]
Expanding each side, we see they are both $-qq'F_2-q^2f_2$.  The remainder of the table is proved in this manner.

For $n=2,4,6$, the identities presented in this table are sufficient to write $c_n$ in terms of $f_2$ and its derivatives.  At $c_8$ we have a potential problem.  There is a nonzero eighth order term appearing in $\fr{1}{2}\left[ (\tr \cT)^2 - \tr \cT^2 \right]$:
\ba
\fr{1}{t^8} \left[u_{1,0}(s_1)^2-u_{0,0}(s_1)u_{1,1}(s_1) \right]\left[u_{1,0}(s_2)^2-u_{0,0}(s_2)u_{1,1}(s_2) \right].
\ea
This term must be multiplied $F_2(s_1)F_2(s_2)$, and the table does not provide a way to reduce this to an $f_2$-type expression.  However, a direct computation will verify the formula
\bae
\left[u_{1,0}^2-u_{0,0}u_{1,1} \right]F_2 = -\fr{1}{6}f_2 +\fr{1}{3} sf_2' - \fr{1}{12}f_2^{(iii)}. \label{8th id}
\eae
Using this formula, the table of identities, and \eqref{scriptt} and \eqref{scriptt2} we arrive at the expression\footnote{We do not explicitly provide the formula for $c_{10}$, as the result is rather cumbersome.}
\bae
c_8(s_1,s_2) &=& \left[ \fr{149}{6} + s_1^3 + 7s_1s_2^2+ 7s_1^2s_2 + s_2^3 \right]f_2(s_1)f_2(s_2)+ \nonumber \\
& \ & \left[ 15s_1 +  \fr{34}{3}s_2 \right]f_2(s_1)f_2'(s_2) + \left[ 15s_2 +  \fr{34}{3}s_1 \right]f_2(s_2)f_2'(s_1) + \nonumber \\
& \ & \left[ \fr{3}{2}s_1^2 +  \fr{13}{3}s_1s_2 +  \fr{3}{2}s_2^2 \right]f_2'(s_1)f_2'(s_2) + 3 [ f''_2(s_1)f'_2(s_2) + f''_2(s_2)f'_2(s_1) ] \nonumber \\
& \ & \fr{1}{2}(s_1+s_2) f''_2(s_1)f''_2(s_2) + \fr{1}{24}f_2^{(iii)}(s_1)f_2^{(iii)}(s_2)	\label{c8}.
\eae
Performing the same sort of integration as in \eqref{C4}, \eqref{C6} results in \eqref{1stC8}.

We can find $t^{-10}$ order terms in trace formulas of the form $\tr \cT^j (\tr \cT)^k$ for $j+k \leq 5$ in the expansion \eqref{det}.  Remarkably, when the computation of tenth order terms in \eqref{det} is carried out it is found that all terms cancel except those contained in
\ba
-\tr \cT + \fr{1}{2}\left[ (\tr \cT)^2 - \tr \cT^2 \right].
\ea
The terms needed from $-\tr \cT$ can be dealt with using the table, as before.  The same is not true of $\fr{1}{2}[ (\tr \cT)^2 - \tr \cT^2]$.  There is a great amount of cancellation within this term, so much so that by introducing the notation
\ba
A(s) &=& u_{1,0}(s)^2-u_{0,0}(s)u_{1,1}(s)\\
B(s) &=& -2u_{2,0}(s)^2+u_{1,1}(s)u_{2,0}(s)-u_{1,0}(s) u_{2,1}(s) 2u_{0,0}(s) u_{2,2}(s) \\
   & \ & + 3u_{1,0}(s)u_{3,0}(s) - 3u_{0,0}(s) u_{3,1} (s)\\
C(s) &=& u_{1,0}(s)u_{2,0}(s)-u_{0,0}(s) u_{2,1}(s)
\ea
we find the tenth order term in $\fr{1}{2}[ (\tr \cT)^2 - \tr \cT^2]$ is given by
\ba
2A(s_1)F_2(s_1)B(s_2)F_2(s_2) + 2A(s_2)F_2(s_2)B(s_1)F_2(s_1) + 2C(s_1)F_2(s_1)C(s_2)F_2(s_2).
\ea
Analogous to the above table and \eqref{8th id}, we may verify the following identities in the same fashion:
\ba
A(s)F_2(s) &=&  - \fr{1}{6}f_2(s) + \fr{1}{3}sf'_2(s) - \fr{1}{12}f_2^{(3)}(s) \\
B(s)F_2(s) &=&  - \fr{1}{3}sf_2(s) + \fr{2}{3}s^2f'_2(s) - \fr{1}{6}f_2^{(3)}(s) \\
C(s)F_2(s) &=&  \fr{1}{12}f'_2(s) + \fr{1}{6}sf_2''(s) - \fr{1}{24}f_2^{(4)}(s).
\ea
Using these formulas to integrate terms in $\fr{1}{2}[ (\tr \cT)^2 - \tr \cT^2]$, and the table of identities to do the same for terms in $\tr \cT$, we find \eqref{1stC10} for $C_{10}$.

\bc \textsc{Appendix 2: Numerical Comparisons}\ec
\par
Folkmar Bornemann has generously provided us with a table of values for the Airy$_2$ covariance obtained by directly computing appropriate Fredholm determinants.  See \cite{Bo1, Bo2} for details on the numerical methods used.  We denote these values $\text{cov}_B(t)$ and consider the information for $5 \leq t \leq 25$. Using the following high precision computations courtesy of Pr\"{a}hofer,
\ba
\mu_1 &=& -1.771\,086\,807\,411\,601\,626\ldots \\
\mu_2 &=& 3.949\,943\,272\,220\,377\,513\ldots \\
\mu_3 &=& -9.711\,844\,753\,027\,647\,354\ldots \\
\mu_4 &=& 26.025\,435\,426\,839\,994\,565\ldots
\ea
together with \eqref{1stC4}--\eqref{1stC10}, we obtain the expressions
\ba
C_4 &=& 2\mu_1 = -3.542\,173\,614\,823\,203\,252\ldots \\
C_6 &=& 2\mu_2 + \fr{10}{3}\mu_1^2 = 18.355\,714 \,809\,065\,487 \ldots \\
C_8 &=& 2\mu_3 + 14\mu_2\mu_1 + \fr{13}{2} = -110.863 \,383 \,378 \,407 \,421 \ldots \\
C_{10} &=& 2\mu_4 +24\mu_3\mu_1 + \fr{126}{5}\mu_2^2 + 116\mu_1 = 652.588 \,990 \,733 \,866 \,004\cdots.
\ea
The following table collects values of $\text{cov}_B(t)$, the approximations $\text{cov}_{2,n}(t)$, and the error between the two values measured by $\text{cov}_B(t)-\text{cov}_{2,n}(t)$.
\ba \footnotesize
\begin{array}{|c|ccccccc|}
\hline t & \text{cov}_B(t) & \text{cov}_{2,6}(t) & \text{Error} & \text{cov}_{2,8}(t) & \text{Error} & \text{cov}_{2,10}(t) & \text{Error} \\ \hline
5 & .03527955721 & .03550728796 & -2 \times 10^{-4} & .03522347770 & 6 \times 10^{-5} & .03529030281 & -1 \times 10^{-5} \\
10 & .00966309240 & .00966413835 & -1 \times 10^{-6} & .00966302972 & 6 \times 10^{-8} & .00966309498 & -3 \times 10^{-9} \\
15 & .004376044913 & .00437608706 & -4 \times 10^{-8} & .00437604380 & 1 \times 10^{-9} & .00437604493 & -2 \times 10^{-11} \\
20 & .002478143955 & .00247814822 & -4\times 10^{-9} & .00247814389 & 6 \times 10^{-11} & .00247814396 & -1 \times 10^{-12} \\
25 & .001591006500 & .00159100722 & -7 \times 10^{-10} & .00159100649 & 6 \times 10^{-12} & .00159100650 & -8 \times 10^{-13} \\ \hline
\end{array}
\ea

\end{document}